\documentclass[11pt,twoside]{article}
\usepackage{amsfonts}
\usepackage{amsmath}
{\markboth{\sl  C. Abdelkefi and M.
Rachdi}{\sl
 Weighted norm inequalities}
\newtheorem{thm}{Theorem}[section]
\newtheorem{rmk}{Remark}[section]
\newtheorem{ex}{Example}[section]

\newtheorem{cor}{Corollary}[section]

\numberwithin {equation}{section}
\newenvironment{pf}{\textbf{Proof.}} {\hfill $\Box$}

\begin{document}
\title{\LARGE\bf Classes of operators on weighted function spaces in Dunkl analysis }
\author{Chokri Abdelkefi and Mongi Rachdi
 \footnote{\small This work was
completed with the support of the DGRST research project LR11ES11
and the program CMCU 10G / 1503.}\\ \small  Department of
Mathematics, Preparatory
Institute of Engineer Studies of Tunis  \\ \small 1089 Monfleury Tunis, University of Tunis, Tunisia\\
  \small E-mail : chokri.abdelkefi@ipeit.rnu.tn \\  \small E-mail : rachdi.mongi@ymail.com}%
\date{}
\maketitle
\begin{abstract}
 For indices $p$ and $q$, $1< p\leq q<+\infty$ and a linear
operator $\mathcal{L}$ satisfying some weak-type boundedness
conditions on suitable function spaces, we give in the Dunkl setting
sufficient conditions on non-negative pairs of weight functions to
obtain weighted norm inequalities for the operator $\mathcal{L}$. We
apply our results to obtain weighted $L^p\rightarrow L^q$
boundedness of the Riesz potentials and of the related fractional
maximal operators for the Dunkl transform. Finally, we prove a
weighted generalized Sobolev inequality.
\end{abstract}
\section{Introduction}
\par
 Dunkl theory generalizes classical Fourier analysis on
$\mathbb{R}^d$. The Dunkl operators introduced by C.F. Dunkl in [9]
are differential-difference operators $T_i, 1 \leq i \leq d$
associated to an arbitrary finite reflection group $W$ on
$\mathbb{R}^d$. These operators attached with a root system $R$ and
a non negative multiplicity function $k$, can be considered as
perturbations of the usual partial derivatives by reflection parts.
They provide a useful tool in the study of special functions with
root systems. Moreover, the commutative algebra generated by these
operators has been used in the study of certain exactly solvable
models of quantum mechanics, namely the Calogero-Moser-Sutherland
models, which deal with systems of identical particles in a
one-dimensional space (see [21]). Dunkl theory was further developed
by several mathematicians (see [8, 11, 15, 17]), later was applied
and generalized in different ways by many authors (see [1, 2, 3, 4,
18]). The Dunkl kernel $E_k$ has been introduced by C.F. Dunkl in
[10]. We use the Dunkl kernel to define the Dunkl transform
$\mathcal{F}_k$, which enjoys properties similar to those of the
classical Fourier transform $\mathcal{F}$. If the parameter
$k\equiv0$ then $\mathcal{F}_k$ becomes $\mathcal{F}$ and the $T_i,
1 \leq i \leq d$ reduce to the corresponding partial derivatives
$\frac{\partial}{\partial x_i}, 1 \leq i \leq d$. Therefore Dunkl
analysis can be viewed as a generalization of classical Fourier
analysis (see next section, Remark 2.1). The classical Fourier
transform behaves well with the translation operator $f \mapsto
f(.-y)$, which leaves the Lebesgue measure on $\mathbb{R}^d$
invariant. However, the weighted measure associated to the Dunkl
operators is no longer invariant under the usual translation. One
ends up with the Dunkl translation operators $\tau_x$, $x\in
\mathbb{R}^d$, introduced by K. Trim\`eche in [20] on the space of
infinitely differentiable functions on $\mathbb{R}^d$. An explicit
formula for the Dunkl translation $\tau_x(f)$ of a radial function
$f$ is known. In particular, the boundedness of $\tau_x$ is
established in this case. As a result one obtains a formula for the
convolution $\ast_k$
(see next section).  \\

In this paper, we consider a linear operator $\mathcal{L}$
satisfying some weak-type boundedness conditions on suitable
function spaces. Our aim is to give for $1<p\leq q<+\infty$,
sufficient conditions on the decreasing rearrangement of
non-negative locally integrable weight functions $u$, $v$ on
$\mathbb{R}^d$, such that $\mathcal{L}$ satisfies the weighted
inequality
\begin{eqnarray}\Big(\int_{\mathbb{R}^{d}}|\mathcal{L}(f)(y)|^{q}u(y)
d\nu_{k}(y)\Big)^{\frac{1}{q}} \leq
c\,\Big(\int_{\mathbb{R}^{d}}|f(x)|^{p}v(x)
d\nu_{k}(x)\Big)^{\frac{1}{p}},\end{eqnarray} for $f\in
L^p_{k,v}(\mathbb{R}^d)$ with $L^p_{k,v}(\mathbb{R}^d)$ is the space
$L^{p}(\mathbb{R}^d, v(x) d\nu_k(x))$ and $\nu_k$ the weighted
measure associated to the Dunkl operators defined by
\begin{eqnarray*}d\nu_k(x):=w_k(x)dx\quad
\mbox{where}\;\;w_k(x) = \prod_{\xi\in R_+} |\langle
\xi,x\rangle|^{2k(\xi)}, \quad x \in \mathbb{R}^d.\end{eqnarray*}
$\langle .,.\rangle$ being the standard Euclidean scalar product on
$\mathbb{R}^d$ and $R_+$ a positive root system (see next section).
As applications of our results, we consider for $0<\alpha<2\gamma+d$
with $\displaystyle\gamma= \sum_{\xi \in R_+}k(\xi),$ the Riesz
potential $I_{\alpha}^{k}$ associated to the Dunkl transform,
defined on the Schwartz space $S(\mathbb{R}^{d})$ by
 \begin{eqnarray}
I_{\alpha}^{k}f(x)=2^{\gamma+\frac{d}{2}-\alpha}\frac{\Gamma(\gamma+\frac{d-\alpha}{2})}{\Gamma(\frac{\alpha}{2})}
\int_{\mathbb{R}^{d}}\frac{\tau_{y}f(x)}{\|y\|^{2\gamma+d-\alpha}}d\nu_{k}(y).
\end{eqnarray}
\\We apply the inequality (1.1) to obtain weighted $L^p\rightarrow L^q$ boundedness of the Riesz potential
$I_\alpha^k$ and of the related fractional maximal operator
$M_{k,\alpha}$ given by
\begin{eqnarray}M_{k,\alpha}f(x)= \displaystyle \sup_{r>0}\frac{1}{m_{k}\,r^{d+2\gamma-\alpha}}
 \int_{\mathbb{R}^{d}}|f(y)|\,\tau_{x}\chi_{B_{r}}(y)d\nu_{k}(y),\quad x\in \mathbb{R}^{d},\end{eqnarray}
where $$m_{k}=\Big(c_k\,
2^{\gamma+\frac{d}{2}}\Gamma(\gamma+\frac{d}{2}+1)\Big)^{\frac{\alpha}{d+2\gamma}-1}\quad\mbox{with}\quad
c_k=\Big(\int_{\mathbb{R}^d} e^{- \frac{\|x\|^2}{2}} w_k
(x)dx\Big)^{-1}.$$ $\chi_{B_{r}}$ denotes the characteristic
function of the ball $B_{r}$ of radius $r$ centered at $0$ (see next
section). These results are generalizations of those obtained by
Heinig for the
classical case in [14]. Finally, we prove a weighted generalized Sobolev inequality.\\

The contents of this paper are as follows. \\In section 2, we
collect some basic definitions and results about harmonic analysis
associated with Dunkl operators .\\
The section 3 is devoted to the proof of the weighted
$L^p\rightarrow L^q$ boundedness of the linear operator
$\mathcal{L}$ when the weights $u,v$ satisfy some conditions.
\\In section 4, as examples, we apply our results to the Riesz potential operators and the related fractional maximal operators.
We obtain finally a weighted generalized Sobolev inequality.\\

Along this paper, we use $c$ to denote a suitable positive constant
which is not necessarily the same in each occurrence and we write
for $x \in \mathbb{R}^d, \|x\| = \sqrt{\langle x,x\rangle}$.
Furthermore, we denote by

$\bullet\quad \mathcal{E}(\mathbb{R}^d)$ the space of infinitely
differentiable functions on $\mathbb{R}^d$.

$\bullet\quad \mathcal{S}(\mathbb{R}^d)$ the Schwartz space of
functions in $\mathcal{E}( \mathbb{R}^d)$ which are rapidly
decreasing as well as their derivatives.

$\bullet\quad \mathcal{D}(\mathbb{R}^d)$ the subspace of
$\mathcal{E}(\mathbb{R}^d)$ of compactly supported functions.
\section{Preliminaries}
$ $ In this section, we recall some notations and results in Dunkl
theory and we refer for more details to the surveys [16].\\

Let $W$ be a finite reflection group on $\mathbb{R}^{d}$, associated
with a root system $R$. For $\alpha\in R$, we denote by
$\mathbb{H}_\alpha$ the hyperplane orthogonal to $\alpha$. For a
given $\beta\in\mathbb{R}^d\backslash\bigcup_{\alpha\in R}
\mathbb{H}_\alpha$, we fix a positive subsystem $R_+=\{\alpha\in R:
\langle \alpha,\beta\rangle>0\}$. We denote by $k$ a nonnegative
multiplicity function defined on $R$ with the property that $k$ is
$W$-invariant. We associate with $k$ the index
$$\gamma = \sum_{\xi \in R_+} k(\xi),$$
and a weighted measure $\nu_k$ given by
\begin{eqnarray*}d\nu_k(x):=w_k(x)dx\quad
\mbox{ where }\;\;w_k(x) = \prod_{\xi\in R_+} |\langle
\xi,x\rangle|^{2k(\xi)}, \quad x \in \mathbb{R}^d,\end{eqnarray*}

Further, we introduce the Mehta-type constant $c_k$ by
$$c_k = \left(\int_{\mathbb{R}^d} e^{- \frac{\|x\|^2}{2}}
w_k (x)dx\right)^{-1}.$$

For every $1 \leq p \leq + \infty$, we denote respectively by
$L^p_k(\mathbb{R}^d)$, $L^p_{k,u}(\mathbb{R}^d)$,
$L^p_{k,v}(\mathbb{R}^d)$ the spaces $L^{p}(\mathbb{R}^d,
d\nu_k(x)),$ $L^{p}(\mathbb{R}^d, u(x)d\nu_k(x)),$
$L^{p}(\mathbb{R}^d, v(x)d\nu_k(x))$ and we use respectively $\|\
\;\|_{p,k}$\,, $\|\ \;\|_{p,k,u}$\,, $\|\ \;\|_{p,k,v}$
as a shorthand for $\|\ \;\|_{L^p_k( \mathbb{R}^d)}$, $\|\ \;\|_{L^p_{k,u}( \mathbb{R}^d)}$, $\|\ \;\|_{L^p_{k,v}( \mathbb{R}^d)}.$ \\

By using the homogeneity of degree $2\gamma$ of $w_k$, it is shown
in [15] that for a radial function $f$ in $L^1_k ( \mathbb{R}^d)$,
there exists a function $F$ on $[0, + \infty)$ such that $f(x) =
F(\|x\|)$, for all $x \in \mathbb{R}^d$. The function $F$ is
integrable with respect to the measure $r^{2\gamma+d-1}dr$ on $[0, +
\infty)$ and we have
 \begin{eqnarray} \int_{\mathbb{R}^d}  f(x)\,d\nu_k(x)&=&\int^{+\infty}_0
\Big( \int_{S^{d-1}}f(ry)w_k(ry)d\sigma(y)\Big)r^{d-1}dr\nonumber\\
&=&
 \int^{+\infty}_0
\Big( \int_{S^{d-1}}w_k(ry)d\sigma(y)\Big)
F(r)r^{d-1}dr\nonumber\\&= & d_k\int^{+ \infty}_0 F(r)
r^{2\gamma+d-1}dr,
\end{eqnarray}
  where $S^{d-1}$
is the unit sphere on $\mathbb{R}^d$ with the normalized surface
measure $d\sigma$  and \begin{eqnarray}d_k=\int_{S^{d-1}}w_k
(x)d\sigma(x) = \frac{c^{-1}_k}{2^{\gamma +\frac{d}{2} -1}
\Gamma(\gamma + \frac{d}{2})}\;.  \end{eqnarray}

The Dunkl operators $T_j\,,\ \ 1\leq j\leq d\,$, on $\mathbb{R}^d$
associated with the reflection group $W$ and the multiplicity
function $k$ are the first-order differential- difference operators
given by
$$T_jf(x)=\frac{\partial f}{\partial x_j}(x)+\sum_{\alpha\in R_+}k(\alpha)
\alpha_j\,\frac{f(x)-f(\rho_\alpha(x))}{\langle\alpha,x\rangle}\,,\quad
f\in\mathcal{E}(\mathbb{R}^d)\,,\quad x\in\mathbb{R}^d\,,$$ where
$\rho_\alpha$ is the reflection on the hyperplane
$\mathbb{H}_\alpha$ and $\alpha_j=\langle\alpha,e_j\rangle,$
$(e_1,\ldots,e_d)$ being the canonical basis of
$\mathbb{R}^d$.
\begin{rmk}In the case $k\equiv0$, the weighted function $w_k\equiv1$ and the measure $\nu_k$ associated to the
Dunkl operators coincide with the Lebesgue measure. The $T_j$ reduce
to the corresponding partial derivatives. Therefore Dunkl analysis
can be viewed as a generalization of classical Fourier
analysis.\end{rmk}

For $y \in \mathbb{C}^d$, the system
$$\left\{\begin{array}{lll}T_ju(x,y)&=&y_j\,u(x,y),\qquad1\leq j\leq d\,,\\  &&\\
u(0,y)&=&1\,.\end{array}\right.$$ admits a unique analytic solution
on $\mathbb{R}^d$, denoted by $E_k(x,y)$ and called the Dunkl
kernel. This kernel has a unique holomorphic extension to
$\mathbb{C}^d \times \mathbb{C}^d $. We have for all $\lambda\in
\mathbb{C}$ and $z, z'\in \mathbb{C}^d,\;
 E_k(z,z') = E_k(z',z)$,  $E_k(\lambda z,z') = E_k(z,\lambda z')$ and for $x, y
\in \mathbb{R}^d,\;|E_k(x,iy)| \leq 1$.\\

The Dunkl transform $\mathcal{F}_k$ is defined for $f \in
\mathcal{D}( \mathbb{R}^d)$ by
$$\mathcal{F}_k(f)(x) =c_k\int_{\mathbb{R}^d}f(y) E_k(-ix, y)d\nu_k(y),\quad
x \in \mathbb{R}^d.$$  We list some known properties of this
transform:
\begin{itemize}
\item[i)] The Dunkl transform of a function $f
\in L^1_k( \mathbb{R}^d)$ has the following basic property
\begin{eqnarray*}\| \mathcal{F}_k(f)\|_{\infty,k} \leq
 \|f\|_{ 1,k}\;. \end{eqnarray*}
\item[ii)] The Dunkl transform is an automorphism on the Schwartz space $\mathcal{S}(\mathbb{R}^d)$.
\item[iii)] When both $f$ and $\mathcal{F}_k(f)$ are in $L^1_k( \mathbb{R}^d)$,
 we have the inversion formula \begin{eqnarray*} f(x) =   \int_{\mathbb{R}^d}\mathcal{F}_k(f)(y) E_k( ix, y)d\nu_k(y),\quad
x \in \mathbb{R}^d.\end{eqnarray*}
\item[iv)] (Plancherel's theorem) The Dunkl transform on $\mathcal{S}(\mathbb{R}^d)$
 extends uniquely to an isometric automorphism on
$L^2_k(\mathbb{R}^d)$.
\item[v)] For $f\in\mathcal{S}(\mathbb{R}^d)$ and $1\leq j \leq d $, we have
\begin{eqnarray}\mathcal{F}_{k}(T_{j}f)(\xi)=i\xi_{j}\mathcal{F}_{k}(f)(\xi), \;\;\xi\in\mathbb{R}^{d}.\end{eqnarray}
\end{itemize}

K. Trim\`eche has introduced in [20] the Dunkl translation operators
$\tau_x$, $x\in\mathbb{R}^d$, on $\mathcal{E}( \mathbb{R}^d)\;.$ For
$f\in \mathcal{S}( \mathbb{R}^d)$ and $x, y\in\mathbb{R}^d$, we have
\begin{eqnarray*}\mathcal{F}_k(\tau_x(f))(y)=E_k(i x, y)\mathcal{F}_k(f)(y).\end{eqnarray*}  Notice that for all $x,y\in\mathbb{R}^d$,
$\tau_x(f)(y)=\tau_y(f)(x)$ and for fixed $x\in\mathbb{R}^d$
\begin{eqnarray*}\tau_x \mbox{\; is a continuous linear mapping from \;}
\mathcal{E}( \mathbb{R}^d) \mbox{\;
into\;}\mathcal{E}(\mathbb{R}^d)\,.\end{eqnarray*} As an operator on
$L_k^2(\mathbb{R}^d)$, $\tau_x$ is bounded. A priori it is not at
all clear whether the translation operator can be defined for $L^p$-
functions with $p$ different from 2. However, according to ([18],
Theorem 3.7), the operator $\tau_x$ can be extended to the space of
radial functions $L^p_k(\mathbb{R}^d)^{rad},$ $1 \leq p \leq 2$ and
we have for a function $f$ in $L^p_k(\mathbb{R}^d)^{^{rad}}$,
\begin{eqnarray}\|\tau_x(f)\|_{p,k} \leq \|f\|_{p,k}.\end{eqnarray}

The Dunkl convolution product $\ast_k$ of two functions $f$ and $g$
in $L^2_k(\mathbb{R}^d)$ is given by
\begin{eqnarray*}(f\; \ast_k g)(x) = \int_{\mathbb{R}^d} \tau_x (f)(-y) g(y) d\nu_k(y),\quad
x \in \mathbb{R}^d.\end{eqnarray*} The Dunkl convolution product is
commutative and for $f,\,g \in \mathcal{D}( \mathbb{R}^d)$, we have
\begin{eqnarray}\mathcal{F}_k(f\,\ast_k\, g) =
\mathcal{F}_k(f) \mathcal{F}_k(g).\end{eqnarray} It was shown in
([18], Theorem 4.1) that when $g$ is a bounded radial function in
$L^1_k( \mathbb{R}^d)$, then
\begin{eqnarray*}(f\; \ast_k g)(x) = \int_{\mathbb{R}^d}  f(y) \tau_x (g)(-y) d\nu_k(y),\quad
x \in \mathbb{R}^d , \end{eqnarray*} initially defined on the
intersection of $L^1_k(\mathbb{R}^d)$ and $L^2_k(\mathbb{R}^d)$
extends to $L^p_k(\mathbb{R}^d)$, $1\leq p\leq +\infty$ as a bounded
operator. In particular, \begin{eqnarray*}\|f \ast_k g\|_{p,k} \leq
\|f\|_{p,k} \|g\|_{1,k}.\end{eqnarray*}
\section{Weighted norm inequalities}
In this section, we prove for a linear operator $\mathcal{L}$
satisfying some weak-type boundedness, weighted norm inequalities
with sufficient conditions on non-negative pairs of weight
functions. We denote by $p'$ the conjugate of $p$ for $1<p<+\infty$.
The proof requires a useful well-known facts which we shall now
state in the following remark.
\begin{rmk}$ $\\
1/ (see [6]) (Hardy inequalities) If $\mu$ and $\vartheta$ are
locally integrable
 weight functions on $(0,+\infty)$ and $1<p\leq q<+\infty$, then there is a constant $c>0$ such that for
all non-negative Lebesgue measurable function $f$ on $(0,+\infty)$,
the inequality
\begin{eqnarray}\Big(\int_{0}^{+\infty}\Big[\int_{0}^{t}f(s)ds\Big]^{q}\mu(t)dt\Big)^{\frac{1}{q}}\leq
 c\, \Big(\int_{0}^{+\infty}(f(t))^{p}\vartheta(t)dt\Big)^{\frac{1}{p}}\end{eqnarray}
 is satisfied if and only if
\begin{eqnarray}\displaystyle\sup_{s>0}\Big(\int_{s}^{+\infty}\mu(t)dt\Big)^{\frac{1}{q}}
\Big(\int_{0}^{s}(\vartheta(t))^{1-p'}dt\Big)^{\frac{1}{p'}}<+
\infty.
\end{eqnarray}
Similarly for the dual operator,
\begin{eqnarray}\Big(\int_{0}^{+\infty}\Big[\int_{t}^{+\infty}f(s)ds\Big]^{q}\mu(t)dt\Big)^{\frac{1}{q}}\leq
 c\Big(\int_{0}^{+\infty}(f(t))^{p}\vartheta(t)dt\Big)^{\frac{1}{p}}
 \end{eqnarray}
 is satisfied if and only if
\begin{eqnarray}\displaystyle\sup_{s>0}\Big(\int_{0}^{s}\mu(t)dt\Big)^{\frac{1}{q}}
\Big(\int_{s}^{+\infty}(\vartheta(t))^{1-p'}dt\Big)^{\frac{1}{p'}}<+
\infty.
\end{eqnarray}
2/ Let $f$ be a complex-valued $\nu_k$-measurable function on
$\mathbb{R}^{d}$. The distribution function $D_f$ of $f$ is defined
for all $s\geq0$ by
$$D_{f}(s)=\nu_k(\{x\in\mathbb{R}^{d}\,:\; |f(x)|>s\}).$$ The decreasing
rearrangement of $f$ is the function $f^*$ given for all $t\geq0$ by
$$f^{*}(t)=inf\{s\geq0 \,:\; D_{f}(s)\leq t\}.$$
We list some known results:
\begin{itemize}
\item[$\bullet$] Let $f\in
L^p_{k}(\mathbb{R}^d)$ and $1\leq p<+\infty$, then
\begin{eqnarray*}\int_{\mathbb{R}^{d}}|f(x)|^{p}
d\nu_{k}(x)=p\int_{0}^{+\infty}s^{p-1}D_{f}(s)ds=\int_{0}^{+\infty}(f^{*}(t))^{p}dt.\end{eqnarray*}
\item[$\bullet$] (see [12]) (Hardy-Littlewood rearrangement inequality)\\
Let $f$ and $\upsilon$ be non negative $\nu_{k}$-measurable
functions on $\mathbb{R}^{d}$, then
\begin{eqnarray}\int_{\mathbb{R}^{d}}f(x)\upsilon(x)d\nu_{k}(x)\leq\int_{0}^{+\infty}f^{*}(t)\upsilon^{*}(t)dt
\end{eqnarray}
and
\begin{eqnarray}\int_{0}^{+\infty}f^{*}(t)\frac{1}
{(\frac{1}{\upsilon})^{*}(t)}dt\leq\int_{\mathbb{R}^{d}}f(x)\upsilon(x)d\nu_{k}(x).
\end{eqnarray}
\item[$\bullet$] (see A. P. Calder\'{o}n [7]) Let $1\leq p_{1}< p_{2}<\infty
$ and $1 \leq q_{1}<q_{2}<+ \infty $. A sublinear operator
$\mathcal{L}$ satisfies the weak-type hypotheses $(p_{1},q_{1})$ and
$(p_{2},q_{2})$  if and only if
\begin{eqnarray}
  (\mathcal{L}f)^{\ast}(t)&\leq& c\,\Big(t^{-\frac{1}{q_{1}}}\int_{0}^{ t^{\frac{\lambda_{1}}{\lambda_{2}}}}s^{\frac{1}{p_{1}}-1}f^{\ast}(s)ds+
   t^{-\frac{1}{q_{2}}}\int_{t^{\frac{\lambda_{1}}{\lambda_{2}}}}^{+\infty}s^{\frac{1}{p_{2}}-1}f^{\ast}(s)ds
   \Big),\nonumber\\ &&
\end{eqnarray} where $\lambda_{1} =\frac{1}{q_{1}}-\frac{1}{q_{2}}$
and $\lambda_{2}=\frac{1}{p_{1}}-\frac{1}{p_{2}}$.
\end{itemize}
\end{rmk}
\begin{ex} Let $\delta<0$ and $\beta>0$. Take $u(x)=\|x\|^{\delta}$, $v(x)=\|x\|^{\beta}$,
 and $f(x)=\chi_{(0,r)}(\|x\|)$, $x\in\mathbb{R}^{d}$, $r>0$. Then using
(2.1) and (2.2), we
have for $s\geq0$, \begin{eqnarray*}D_{u}(s)&=&\nu_{k}\Big(\{x\in\mathbb{R}^{d}\,:\;\|x\|^{\delta}>s\}\Big)\\
&=&\nu_{k}\Big(B(0,s^{\frac{1}{\delta}})\Big)=
\frac{d_{k}}{2\gamma+d}\;s^{\frac{2\gamma+d}{\delta}},\end{eqnarray*}
\begin{eqnarray*}D_{\frac{1}{v}}(s)&=&\nu_{k}\Big(\{x\in\mathbb{R}^{d}\,:\;\|x\|^{-\beta}>s\}\Big)\\
&=&\nu_{k}\Big(B(0,s^{-\frac{1}{\beta}})\Big)=\frac{d_{k}}{2\gamma+d}\;s^{-\frac{2\gamma+d}{\beta}},\end{eqnarray*}
and \begin{eqnarray*}D_{f}(s)&=&\nu_{k}\Big(\{x\in\mathbb{R}^{d}\,:\;\chi_{(0,r)}(\|x\|)>s\}\Big)\\
&=& \nu_{k}\Big(B(0,1)\Big)\,r^{2\gamma+d}\,\chi_{(0,1)}(s)\\
&=&\frac{d_k}{2\gamma+d}\,r^{2\gamma+d}\,\chi_{(0,1)}(s).\end{eqnarray*}
Note that from (2.2), $\displaystyle
d_k=\frac{c^{-1}_k}{2^{\gamma+\frac{d}{2} -1} \Gamma(\gamma +
\frac{d}{2})},$ this yields
$$\frac{d_k}{2\gamma+d}=  \frac{c^{-1}_k}{2^{\gamma +\frac{d}{2}}
\Gamma(\gamma + \frac{d}{2}+1)}.$$ This gives for $t\geq0$,
\begin{eqnarray*}u^{*}(t)=inf\{s\geq 0\,:\; D_{u}(s)\leq t\}=
\Big(\frac{2\gamma+d}{d_{k}}\Big)^{\frac{\delta}{2\gamma+d}}\;t^{\frac{\delta}{2\gamma+d}},\end{eqnarray*}
\begin{eqnarray*}(\frac{1}{v})^{*}(t)=inf\{s\geq 0\,:\; D_{\frac{1}{v}}(s)\leq t\}
=\Big(\frac{2\gamma+d}{d_{k}}\Big)^{-\frac{\beta}{2\gamma+d}}\;t^{-\frac{\beta}{2\gamma+d}},\end{eqnarray*}
and \begin{eqnarray*}f^*(t)=\chi_{(0,R)}(t)\;\;\mbox{where}\;\;
R=\frac{d_k}{2\gamma+d}\,r^{2\gamma+d}.\end{eqnarray*} Hence, using
(2.1) and (2.2) again, we obtain for $-(2\gamma+d)<\delta$
\begin{eqnarray*}\int_{\mathbb{R}^{d}}
f(x)u(x)d\nu_{k}(x)&=& \frac{d_k
}{\delta+2\gamma+d}\;r^{\delta+2\gamma+d}\\&=&\int_{0}^{+\infty}f^{*}(t)u^{*}(t)dt,\end{eqnarray*}
and
\begin{eqnarray*}\int_{0}^{+\infty}f^{*}(t)\frac{1}
{(\frac{1}{v})^{*}(t)}dt&=& \frac{d_k
}{\beta+2\gamma+d}\;r^{\beta+2\gamma+d}\\&=&\int_{\mathbb{R}^{d}}
f(x)v(x)d\nu_{k}(x),\end{eqnarray*} giving equalities for (3.5) and
(3.6) in these cases.
\end{ex}$ $\\

In the following theorem, we prove weighted norm inequalities for a
linear operator.
\begin{thm} Let $u$,$v$ be non negative $\nu_{k}$-locally integrable weight functions on
$\mathbb{R}^{d}$ and $1\leq p_{1}<p_{2}<+\infty$, $ 1\leq
q_{1}<q_{2}<+ \infty $. Suppose $\mathcal{L}$ is a linear operator
defined on $\mathcal{S}(\mathbb{R}^{d})$ such that $\mathcal{L}$
 is simultaneously of weak-type $(p_{i},q_{i})$, $i=1,2$,
 then for $1< p \leq q <+ \infty $, $\mathcal{L}$ can be extended to a bounded operator from $L_{k,v}^{p}(\mathbb{R}^{d})$
  to $L_{k,u}^{q}(\mathbb{R}^{d})$ and the inequality \begin{eqnarray*}\|\mathcal{L}f\|_{q,k,u}\leq c\,\|f\|_{p,k,v}\end{eqnarray*}
   holds with the following conditions on $u$ and
  $v$:
\begin{eqnarray}
  \sup_{s>0}\Big(\int_{s^{\frac{1}{\lambda_{1}}}}^{+\infty}u^{\ast}(t)t^{-\frac{q}{q_{1}}}dt\Big)^{\frac{1}{q}}
  \Big(\int_{0}^{s^{\frac{1}{\lambda}_{2}}}[(\frac{1}{v})^{\ast}(t)]^{(p'-1)}t^{p'(\frac{1}{p_{1}}-1)}dt\Big)^{\frac{1}{p'}}<+\infty
  \end{eqnarray} and \begin{eqnarray}
 \sup_{s>0}\Big(\int_{0}^{s^{\frac{1}{\lambda_{1}}}}u^{\ast}(t)t^{-\frac{q}{q_{2}}}dt\Big)^{\frac{1}{q}}
 \Big(\int_{s^{\frac{1}{\lambda_{2}}}}^{+\infty}[(\frac{1}{v})^{\ast}(t)]^{(p'-1)}t^{p'(\frac{1}{p_{2}}-1)}dt\Big)^{\frac{1}{p'}}<+\infty,
\end{eqnarray} where $\lambda_{1} = \frac{1}{q_{1}}-\frac{1}{q_{2}}$ and
$\lambda_{2}=\frac{1}{p_{1}}-\frac{1}{p_{2}}$.
 \end{thm}
\begin{pf} Let $f\in\mathcal{S}(\mathbb{R}^{d})$. Using (3.7) and applying Minkowski's inequality, we have\\
$\displaystyle\Big(\int_{0}^{+\infty}[(\mathcal{L}f)^{\ast}(t)]^{q}u^{\ast}(t)dt\Big)^{\frac{1}{q}}$
\begin{eqnarray*}
 &\leq & c\,
   \Big[\int_{0}^{+\infty}u^{\ast}(t)t^{-\frac{q}{q_{1}}}
  \Big(\int_{0}^{t^{\frac{\lambda_{1}}{\lambda_{2}}}}s^{\frac{1}{p_{1}}-1}f^{\ast}(s)ds\Big)^{q}dt\Big]^{\frac{1}{q}} \\
   && +\,c\,\Big[\int_{0}^{+\infty}u^{\ast}(t)t^{-\frac{q}{q_{2}}}
  \Big(\int_{t^{\frac{\lambda_{1}}{\lambda_{2}}}}^{+\infty}s^{\frac{1}{p_{2}}-1}f^{\ast}(s)ds\Big)^{q}dt\Big]^{\frac{1}{q}}.
\end{eqnarray*}
By means of change of variable in the right side, we obtain\\
$\displaystyle\Big(\int_{0}^{+\infty}[(\mathcal{L}f)^{\ast}(t)]^{q}u^{\ast}(t)dt\Big)^{\frac{1}{q}}$
\begin{eqnarray}
   &\leq & c\,
 \Big[\int_{0}^{+\infty}u^{\ast}(t^{\frac{\lambda_{2}}{\lambda_{1}}})t^{\frac{\lambda_{2}}{\lambda_{1}}(1-\frac{q}{q_{1}})-1}
  [\int_{0}^{t}s^{\frac{1}{p_{1}}-1}f^{\ast}(s)ds]^{q}dt\Big]^{\frac{1}{q}}\nonumber\\
  &&+\,
  c\,\Big[\int_{0}^{+\infty}u^{\ast}(t^{\frac{\lambda_{2}}{\lambda_{1}}})t^{\frac{\lambda_{2}}{\lambda_{1}}(1-\frac{q}{q_{2}})-1}
  [\int_{t}^{+\infty}s^{\frac{1}{p_{2}}-1}f^{\ast}(s)ds]^{q}dt\Big]^{\frac{1}{q}}\nonumber \\
 & =& I_{1}+I_{2}.
 \end{eqnarray}
Applying (3.1) and (3.2) for $I_1$, we can assert that
 \begin{eqnarray} I_{1}\leq \Big(\int_{0}^{+\infty}[(\frac{1}{v})^{\ast}(t)]^{-1}[f^{\ast}(t)]^{p}dt\Big)^{\frac{1}{p}}\end{eqnarray}
 if and only if
 \begin{eqnarray*}
   \sup_{s>0}\Big(\int_{s}^{+\infty}u^{\ast}(t^{\frac{\lambda_{2}}{\lambda_{1}}})t^{\frac{\lambda_{2}}{\lambda_{1}}(1-\frac{q}{q_{1}})-1}dt\Big)^{\frac{1}{q}}
    \Big(\int_{0}^{s}[(\frac{1}{v})^{\ast}(t)]^{(p'-1)}t^{p'(\frac{1}{p_{1}}-1)}dt\Big)^{\frac{1}{p'}}\leq
    +\infty.
 \end{eqnarray*}
Then if we replace $s$ by $s^{\frac{1}{\lambda_{2}}}$ in this
condition, it's easy to see that if we use a change of variable in
the first integral of the expression, we obtain (3.8).  \\
 Similarly by applying (3.3) and (3.4) for $I_2$, we get
\begin{eqnarray}
    I_{2}\leq \Big(\int_{0}^{+\infty}[(\frac{1}{v})^{\ast}(t)]^{-1}[f^{\ast}(t)]^{p}dt\Big)^{\frac{1}{p}}
 \end{eqnarray}
 if and only if
 \begin{eqnarray*}
    \sup_{s>0}\Big(\int_{0}^{s}u^{\ast}(t^{\frac{\lambda_{2}}{\lambda_{1}}})t^{\frac{\lambda_{2}}{\lambda_{1}}(1-\frac{q}{q_{2}})-1}dt\Big)^{\frac{1}{q}}
    \Big(\int_{s}^{+\infty}[(\frac{1}{v})^{\ast}(t)]^{(p'-1)}t^{p'(\frac{1}{p_{2}}-1)}dt\Big)^{\frac{1}{p'}}\leq
    +\infty,
 \end{eqnarray*}
 which is equivalent to (3.9).\\
 Combining (3.10), (3.11) and (3.12), it yields
 \begin{eqnarray}
 \displaystyle\Big(\int_{0}^{+\infty}[(\mathcal{L}f)^{\ast}(t)]^{q}u^{\ast}(t)dt\Big)^{\frac{1}{q}}
  \leq
  c\,\Big(\int_{0}^{+\infty}[(\frac{1}{v})^{\ast}(t)]^{-1}[f^{\star}(t)]^{p}dt\Big)^{\frac{1}{p}}.
  \end{eqnarray}
 Using (3.5) on the left side and (3.6) on the right side of (3.13), we
 obtain by density of $\mathcal{S}(\mathbb{R}^{d})$ in
 $L_{k,v}^{p}(\mathbb{R}^{d})$, $1\leq p<+\infty$,
 \begin{eqnarray*}
 \Big(\int_{\mathbb{R}^{d}}[(\mathcal{L}f)(x)]^{q}u(x)d\nu_{k}(x)\Big)^{\frac{1}{q}}
 \leq c\,
 \Big(\int_{\mathbb{R}^{d}}[f(x)]^{p}v(x)d\nu_{k}(x)\Big)^{\frac{1}{p}}.
  \end{eqnarray*}
This completes the proof.
\end{pf}
\section{Applications}$ $
Our first application is for the Riesz Potential $I_{\alpha}^{k}$.
It was shown in [13] that the Riesz potential $I_{\alpha}^{k}$
satisfies the following Hardy-Littlewood-Sobolev theorem.
\begin{thm} Let $0<\alpha<2\gamma+d$. Then
\begin{itemize}
  \item[i)] for $f\in L_{k}^{r}(\mathbb{R}^{d})$, $1<
r<\frac{2\gamma+d}{\alpha}$, the mapping $f\longrightarrow
I_{\alpha}^{k}f $ is of strong-type
  $\Big(r\,,\,\ell=\frac{1}{\frac{1}{r}-\frac{\alpha}{2\gamma+d}}\Big)$
  and one has
  $$\|I_{\alpha}^{k}f\|_{\ell,k}\leq c \,\|f\|_{r,k}.$$
\item[ii)] for $f\in L_{k}^{1}(\mathbb{R}^{d})$, the mapping $f\longrightarrow I_{\alpha}^{k}f $ is of weak-type
$(1,\frac{1}{1-\frac{\alpha}{2\gamma+d}})$ and one has for any
$\lambda>0,$
$$\int_{\{x \in \mathbb{R}^{d}:\,|I_{\alpha}^{k}f(x)|>\lambda\}}d\nu_{k}(x) \leq
c\,\Big(\frac{\|f\|_{1,k}}{\lambda}\Big)^{\frac{1}{1-\frac{\alpha}{2\gamma+d}}}.$$
\end{itemize}
\end{thm}
\begin{rmk} The boundedness of Riesz potentials can be used to establish the
boundedness properties of the fractional maximal operator given by
(1.3). This follows from the fact that
$$M_{k,\alpha}f(x)\leq c\,I_{\alpha}^{k}(|f|)(x),\;x\in \mathbb{R}^{d}.$$  Hence
we deduce for $M_{k,\alpha}$, the same results obtained in Theorem
4.1.
\end{rmk}
Now, applying Theorem 3.1, we obtain the following results.
\begin{thm} Let $0<\alpha<2\gamma+d$, $1<r<\frac{2\gamma+d}{\alpha}$ and $u$,$v$ be non negative $\nu_{k}$-locally integrable weight functions on
$\mathbb{R}^{d}$. Then for $1< p \leq q <+ \infty $,
$I_{\alpha}^{k}$ can be extended to a bounded operator from
$L_{k,v}^{p}(\mathbb{R}^{d})$
  to $L_{k,u}^{q}(\mathbb{R}^{d})$ and the inequality \begin{eqnarray*}\|I_{\alpha}^{k}f\|_{q,k,u}\leq c\,\|f\|_{p,k,v}\end{eqnarray*}
   holds with the following conditions on $u$ and
  $v$:
\begin{eqnarray}
  \sup_{s>0}\Big(\int_{s}^{+\infty}u^{\ast}(t)t^{-q(1-\frac{\alpha}{2\gamma+d})}dt
  \Big)^{\frac{1}{q}} \Big(\int_{0}^{s}[(\frac{1}{v})^{\ast}(t)]^{(p'-1)}dt\Big)^{\frac{1}{p'}}<+\infty
  \end{eqnarray}
   and
   \begin{eqnarray}
   \sup_{s>0} \Big(\int_{0}^{s}u^{\ast}(t)t^{-q(\frac{1}{r}-\frac{\alpha}{2\gamma+d})}dt\Big)^{\frac{1}{q}}
\Big(\int_{s}^{+\infty}[(\frac{1}{v})^{\ast}(t)^{(p'-1)}t^{p'(\frac{1}{r}-1)}dt\Big)^{\frac{1}{p'}}
<+\infty.
\end{eqnarray}
\end{thm}
\begin{pf}
From Theorem 4.1, $I_{\alpha}^{k}$ is of weak-type
$(p_1,q_1)=(1,\frac{1}{1-\frac{\alpha}{2\gamma+d}})$ and strong-type
$(p_2,q_2)=(r,\frac{1}{\frac{1}{r}-\frac{\alpha}{2\gamma+d}})$, then
the result follows immediately from Theorem 3.1 with $\lambda_{1}
=\lambda_{2}=1-\frac{1}{r}$.
\end{pf}
\begin{rmk}Note that if $u=v\equiv1$, the boundedness conditions
(4.1) and (4.2) are valid if and only if
$\;1<p<r<\frac{2\gamma+d}{\alpha}$ and
$\frac{1}{q}=\frac{1}{p}-\frac{\alpha}{2\gamma+d}$. Therefore,
Theorem 4.2 reduces to Theorem 4.1, i).
\end{rmk}

As consequence of Theorem 4.2 for power weights, we obtain the
result below.
\begin{cor}
Let $0<\alpha<2\gamma+d$ and $1<p<\frac{2\gamma+d}{\alpha}$. For
$\delta,\beta$ such that $\delta < 0 $, $0< \beta =\delta+\alpha p
<(2\gamma+d)(p-1)$ and $f\in L_{k,v}^{p}(\mathbb{R}^{d})$ with
$v=\|.\|^\beta$, we have
\begin{eqnarray*}
\Big(\int_{\mathbb{R}^{d}}|I_{\alpha}^{k}f(x)|^{p}\|x\|^{\delta}d\nu_{k}(x)\Big)^{\frac{1}{p}}
\leq c\,
  \Big(\int_{\mathbb{R}^{d}}|f(x)|^{p}\|x\|^{\beta}d\nu_{k}(x)\Big)^{\frac{1}{p}}
\end{eqnarray*}
\end{cor}
\begin{pf}
From Example 3.1, we have for $\delta<0$ and
$\beta>0$\begin{eqnarray*}
u^{\ast}(t)=\Big(\frac{2\gamma+d}{d_{k}}\Big)^{\frac{\delta}{2\gamma+d}}t^{\frac{\delta}{2\gamma+d}}
\quad\mbox{and}\quad
(\frac{1}{v})^{\ast}(t)=\Big(\frac{2\gamma+d}{d_{k}}\Big)^{-\frac{\beta}{2\gamma+d}}t^{-\frac{\beta}{2\gamma+d}},\end{eqnarray*}
then if we take $p=q=r$ in Theorem 4.2, the boundedness conditions
(4.1) and (4.2) are valid if and only if
 $$ \left\{\begin{array}{lll}0<\beta<(2\gamma+d)(p-1),\\
&&\\\beta=\delta+\alpha p\,.\end{array}\right.$$ Under these
conditions and from Theorem 4.2, we obtain our result.
\end{pf}
\begin{rmk} As in Remark 4.1, the boundedness of Riesz potentials can be used to establish the
boundedness properties of the fractional maximal operator
$M_{k,\alpha}$. Then we get for $M_{k,\alpha}$ the same results
obtained in Theorem 4.2 and Corollary 4.1.
\end{rmk}
Our next application concerns the weighted Generalized Sobolev
inequality. Before, we need some useful results that we state in the
following remark.
\begin{rmk}$ $ \\
1/ (see[20]) In Dunkl setting the Riesz transforms  are the
operators $\mathcal{R}_{j}$, $j=1...d,$ defined on
$L^{2}_{k}(\mathbb{R}^{d})$ by
\begin{eqnarray*}
  \mathcal{R}_{j}(f)(x) &=& 2^{\frac{\ell_{k}-1}{2}}\frac{\Gamma(\frac{\ell_{k}}{2})}{\sqrt{\pi}} \,\lim_{\epsilon\rightarrow 0}\int_{\|y\|>\epsilon}\tau_{x}(f)(-y)\frac{y_{j}}{\|y\|^{p_{k}}}d\nu_{k}(y),x\in\mathbb{R}^{d}
\end{eqnarray*}where
\begin{eqnarray*}
 \ell_{k}=2\gamma+d+1.
\end{eqnarray*}
$\bullet$ The Riesz transform $\mathcal{R}_{j}$ is a multiplier
operator with
\begin{eqnarray}
\mathcal{F}_{k}(\mathcal{R}_{j}(f))(\xi) =
\frac{-i\xi_{j}}{\|\xi\|}\,\mathcal{F}_{k}(f)(\xi),\;1\leq j\leq d,
\quad f\in
  \mathcal{S}(\mathbb{R}^{d}).
\end{eqnarray}
$\bullet$ Let $0<\alpha<2\gamma+d $. The identity
\begin{eqnarray}\mathcal{F}_{k}(I_{\alpha}^{k}f)(x)=\|x\|^{-\alpha}\mathcal{F}_{k}(f)(x)\end{eqnarray}
holds in the sense that\begin{eqnarray*}
\int_{\mathbb{R}^{d}}I_{\alpha}^{k}f(x)g(x)d\nu_{k}(x)=\int_{\mathbb{R}^{d}}\mathcal{F}_{k}(f)(x)\|x\|^{-\alpha}\mathcal{F}_{k}(g)(x)d\nu_{k}(x),\end{eqnarray*}
whenever $f,g\in \mathcal{S}(\mathbb{R}^{d})$.\\\\
2/ (see [5]) The Riesz transform $\mathcal{R}_{j}$,$\;1\leq j\leq
d$, can be extended to a bounded operator from
$L_{k}^{p}(\mathbb{R}^{d})$ into it self for $1<p<+\infty$ and we
have
\begin{eqnarray}\|\mathcal{R}_{j}(f)\|_{p,k}\leq c\,\|f\|_{p,k}.\end{eqnarray}
\end{rmk}
Now, using Remark 4.4 and Theorem 4.2, we obtain the following
results.
\begin{thm}(Weighted generalized Sobolev inequality) Let $u$ be a non-negative $\nu_{k}$-locally integrable function on
$\mathbb{R}^{d}$ and $1<r<2\gamma+d$. Then for $1<p\leq q<+\infty$
such that $p<r$ and $f\in \mathcal{D}(\mathbb{R}^{d})$, the
inequality \begin{eqnarray*}\|f\|_{q,k,u} \leq
c\,\|\nabla_{k}f\|_{p,k},\end{eqnarray*} holds with the following
conditions on $u$:
\begin{eqnarray}
\Big(\int_{s}^{\infty}u^{\ast}(t)t^{-q(1-\frac{1}{2\gamma+d})}dt\Big)^{\frac{1}{q}}
\leq c\,s^{\frac{1}{p}-1}
    \end{eqnarray}
   and
   \begin{eqnarray}
    \Big(\int_{0}^{s}u^{\ast}(t)t^{-q(\frac{1}{r}-\frac{1}{2\gamma+d})}dt\Big)^{\frac{1}{q}}
\leq c\,s^{\frac{1}{p}-\frac{1}{r}},
\end{eqnarray} for all $s>0$. Here $\nabla_{k}f=(T_{1}f,...,T_{d}f)$ and
$\displaystyle|\nabla_{k}f|=\Big(\sum_{j=1}^{d}|T_{j}f|^{2}\Big)^{\frac{1}{2}}$.
\end{thm}
\begin{pf}
For $f\in \mathcal{D}(\mathbb{R}^{d})$, we write
\begin{eqnarray*}
  \mathcal{F}_{k}(f)(\xi)=\frac{1}{\|\xi\|}\sum_{j=1}^{d}\frac{-i\xi_{j}}{\|\xi\|}\,i\xi_{j}\,\mathcal{F}_{k}(f)(\xi),
 \end{eqnarray*}then by (2.3) and (4.3), we get
 \begin{eqnarray*}\mathcal{F}_{k}(f)(\xi)&=&\frac{1}{\|\xi\|}\sum_{j=1}^{d}\frac{-i\xi_{j}}{\|\xi\|}\,\mathcal{F}_{k}(T_{j}f)(\xi)\\
 &=&\frac{1}{\|\xi\|}\sum_{j=1}^{d} \mathcal{F}_{k}(\mathcal{R}_{j}(T_{j}f))(\xi)\\
 &=&\frac{1}{\|\xi\|}\mathcal{F}_{k}\Big(\sum_{j=1}^{d} \mathcal{R}_{j}(T_{j}f)\Big)(\xi).
\end{eqnarray*}
This yields from (4.4) that\begin{eqnarray*}
\mathcal{F}_{k}(f)(\xi)=\mathcal{F}_{k}\Big[I_{1}^{k}\Big(\sum_{j=1}^{d}\mathcal{R}_{j}(T_{j}f)\Big)\Big](\xi),
\end{eqnarray*} which gives the following identity,
\begin{eqnarray*}
f=I_{1}^{k}\Big(\sum_{j=1}^{d}\mathcal{R}_{j}(T_{j}f)\Big).
\end{eqnarray*}
Now, observe that the conditions (4.6) and (4.7) are equivalent to
(4.1) and (4.2) with $v\equiv1$ and $\alpha =1$, then using Theorem
4.2, we obtain
\begin{eqnarray*}
  \|f\|_{q,k,u} &=& \|I_{1}^{k}\Big(\sum_{j=1}^{d}\mathcal{R}_{j}(T_{j}f)\Big)\|_{q,k,u}\nonumber \\
   &\leq&c\,\|\mathcal{R}_{j}(\sum_{j=1}^{d}(T_{j}f))\|_{p,k},\end{eqnarray*}
which gives from (4.5) that \begin{eqnarray*}\|f\|_{q,k,u}  &\leq&c\,\|\sum_{j=1}^{d}(T_{j}f)\|_{p,k} \\
   &\leq&c\,\|\nabla_{k}f\|_{p,k}.
\end{eqnarray*}Our result is proved.
\end{pf}
\begin{cor}
Let $1<p<2\gamma+d$ and $1<p\leq q<+\infty$. Then for $\delta < 0$
such that $\delta=q\,[(2\gamma+d)(\frac{1}{p}-\frac{1}{q})-1]$, we
have for $f\in \mathcal{D}(\mathbb{R}^{d})$
\begin{eqnarray*}
\Big(\int_{\mathbb{R}^{d}}|f(x)|^{q}\|x\|^{\delta}d\nu_{k}(x)\Big)^{\frac{1}{q}}
\leq c\,\|\nabla_{k}f\|_{p,k}.
\end{eqnarray*}
\end{cor}
\begin{pf}
For $\delta<0$, if we take $u(x)=\|x\|^\delta$, $x\in \mathbb{R}^{d}
$ in Theorem 4.2, the boundedness conditions (4.6) and (4.7) are
valid if and only if
\begin{eqnarray*}\delta=q\,\Big[(2\gamma+d)\Big(\frac{1}{p}-\frac{1}{q}\Big)-1\Big].\end{eqnarray*} Under this
condition and from Theorem 4.3, we obtain our result.
\end{pf}
\begin{rmk}
The case $u\equiv1$, $1<p<2\gamma+d$ and
$\frac{1}{q}=\frac{1}{p}-\frac{1}{2\gamma+d}$ was obtained in [5]
and gives the generalized Sobolev inequality
\begin{eqnarray*}
\|f\|_{q,k} \leq c\,\|\nabla_{k}f\|_{p,k}.
\end{eqnarray*}
\end{rmk}

\end{document}